%% file: compactmin7.tex
\documentclass[11pt]{smfart}

\input{paquets.tex}
\input{lemme-proposition-theoreme.tex}
\input{numerotation}

\input{macros}

\numberwithin{equation}{subsubsection}
\renewcommand{\theequation}{\thesubsection.\Alph{equation}}

\newcommand{\isolong}{\buildrel{\sim}\over\longrightarrow}

\newcommand{\gsp}{\mathrm{GSp}_{2g}}
\newcommand{\ag}{\mathcal{A}_{g}}
\newcommand{\agn}{\mathcal{A}_{g,n}}
\newcommand{\agK}{\mathcal{A}_{g ,0}}
\newcommand{\agKn}{\mathcal{A}_{g ,n,0}}

\newcommand{\agb}{\overline{\mathcal{A}}_{g}}

\newcommand{\agnb}{\overline{\mathcal{A}}_{g,n}}
\newcommand{\agnm}{\mathcal{A}_{g,n}^*}
\newcommand{\agKb}{\overline{\mathcal{A}}_{g,0}}

\newcommand{\agKnb}{\overline{\mathcal{A}}_{g,n,0}}
\newcommand{\agKnm}{\mathcal{A}_{g,n,0}^*}

\title{Compactification minimale et mauvaise réduction}
\author{Benoît Stroh}
\date{1 novembre 2008}
\email{benoit.stroh@gmail.com}
\address{Institut Élie Cartan, Université Henri Poincaré,  B.P.~239, F-54506 Vandoeuvre-les-Nancy Cedex, France}

\begin{document}
\maketitle

\begin{abstract} Nous construisons la compactification minimale de certaines variétés modulaires de Siegel en leurs places de mauvaise réduction. Ces variétés paramètrent des schémas abéliens principalement polarisés munis d'une structure de niveau parahorique en un nombre premier~$p$, et d'une structure de niveau auxilliaire~;~elles ont mauvaise réduction en $p$. Nous esquissons également une théorie arithmétique des formes modulaires de Siegel associées à ces variétés.
\end{abstract}

\begin{altabstract} We construct the minimal compactification of some modular Siegel varieties at their bad reduction places. These varieties parametrize principally polarized abelian schemes endowed with a parahoric level structure at a prime number $p$, and with an auxiliary level structure~;~such varieties have bad reduction at $p$. We also sketch the arithmetic theory of the associated Siegel modular forms.
\end{altabstract}


Soient $g\geq 2$ et $1\leq s \leq g$ deux entiers, $p$ un nombre premier, $n\geq 3$ un entier non divisible par $p$ et $\mathcal{D}=\{d_1<d_2<\cdots<d_s\}$ un sous-ensemble ordonné de $\{1,\cdots,g\}$. Notons $\agK$ le champ dont les points à valeurs dans un schéma $S$ paramètrent les schémas abéliens $G$ sur $S$ principalement polarisés de genre $g$  munis d'un drapeau
$$H_1\subset H_2 \subset \cdots \subset H_{s-1} \subset H_s \subset G[p]$$
de groupes totalement isotropes, finis et plats sur $S$, tels que $H_i$ soit de rang $p^{d_i}$ pour $1\leq~i\leq~s$~;~on dit que $G$ est muni d'une structure de niveau parahorique de type $\mathcal{D}$ en~$p$. Désignons par $G_0$ le schéma abélien universel sur $\agK$ et notons $\agKn$ le champ sur 
$$\agK\times_{\Spec(\Z)} \Spec(\Z[1/n])$$
qui paramètre les structures de niveau principales sur $G_0[n]$.
Les champs $\agK$ et $\agKn$ sont algébriques, sont lisses sur $\Spec(\Z[1/pn])$ mais ne sont pas lisses sur $\Spec(\Z_p)$. On dit qu'ils ont \textit{bonne réduction} en les nombres premiers qui ne divisent pas $pn$ et qu'ils ont \textit{mauvaise réduction} en~$p$.
L'avantage de considérer $\agKn$ plutôt que $\agK$ est que d'après~\cite{Jacobi@Serre} et~\cite{GIT@Mumford}, $\agKn$ est un schéma quasi-projectif sur $\Spec(\Z[1/n])$. Ce schéma n'est pas projectif, et on aimerait le compactifier. Une variante d'une construction de Faltings et Chai~\cite[ch.~V]{Deg@FaltingsChai} permet d'obtenir une compactification canonique
$$\left(\agKn\times_{\Spec(\Z[1/n])} \Spec(\Z[1/pn])\right)^*$$ de la restriction de $\agKn$ au schéma $\Spec(\Z[1/pn])$ formé des places de {bonne réduction}. Cette compactification est appelée la \textit{compactification minimale} de $\agKn\times\Spec(\Z[1/pn])$. Avant de donner plus de précisions, introduisons quelques notations.
Désignons encore par $G_0$ le schéma abélien universel sur $\agKn$ et notons $\omega_0$ le faisceau inversible sur $\agKn$ des $g$-formes différentielles invariantes sur $G_0$. Rappelons que pour tout $k\in\Z$, une forme modulaire de poids $k$ pour $\agKn\times \Spec(\Z[1/pn])$ est un élément du groupe
$$\mathrm{H}^0\left( \agKn\times \Spec(\Z[1/pn]) \: , \: \omega_0^k \right) \: .$$
La compactification minimale de Faltings et Chai s'obtient comme schéma projectif associé à l'algèbre graduée des formes modulaires de poids positif, soit
$$\left(\agKn\times \Spec(\Z[1/pn])\right)^* \: = \: \mathrm{Proj} \left( \bigoplus_{k\in\mathbb{N}}   \mathrm{H}^0\left( \agKn\times \Spec(\Z[1/pn]) \: , \: \omega_0^k \right) \right) \: .$$

Dans cet article, nous construisons une compactification $\agKnm$ de $\agKn$ sur $\Spec(\Z[1/n])$ qui étend $(\agKn\times \Spec(\Z[1/pn]))^*$ au-dessus de la place de mauvaise réduction $p$. Ainsi, le schéma $\agKnm$ est projectif sur $\Spec(\Z[1/n])$ et vérifie
$$ \agKnm \times_{\Spec(\Z[1/n])} \Spec(\Z[1/pn]) \: \isolong \: \left( \agK\times_{\Spec(\Z[1/n])} \Spec(\Z[1/pn]) \right)^* \: .$$
On l'appelle \textit{compactification minimale} de $\agKn$. Le schéma $\agKnm$ est relié à l'algèbre graduée des formes modulaires pour $\agKn$. Remarquons que dans ce contexte, la notion de poids réserve quelques surprises. En effet, notons $G_i=G_0/H_i$ le quotient de $G_0$ par le $i$-ème sous-groupe universel sur $\agKn$, et notons $\omega_i$ le faisceau inversible sur $\agKn$ des $g$-formes différentielles invariantes sur $G_i$ pour $1\leq i \leq s$. La chaîne d'isogénies
$$G_0 \: \longrightarrow \: G_1 \: \longrightarrow \: \cdots\: \longrightarrow \: G_s $$
sur $\agKn$ induit une chaîne de morphismes
$$\omega_s \: \longrightarrow \: \omega_{s-1} \: \longrightarrow \: \cdots \: \longrightarrow \omega_0 \: ,$$
dont les restrictions à $\Spec(\Z[1/pn])$ sont des isomorphismes. Par contre, ces morphismes ne sont pas des isomorphismes au-dessus de $\Spec(\Z_p)$.
Il existe donc plusieurs faisceaux automorphes intéressants sur $\agKn$, et il y a lieu de tous les considérer. 
Ce phénomène est à relier à la théorie du modèle local~\cite{Gamma0@DeJong}, qui étudie la {mauvaise réduction} de $\agKn$. Définissons les formes modulaires de poids $(k_0,\cdots,k_s)\in\Z^{s+1}$ pour $\agKn$ comme les éléments du groupe
$$\mathrm{H}^0\left( \agKn \: , \: \bigotimes_{i=0}^s \: \omega_i^{k_i} \right)\: .$$
On obtient le résultat suivant, qui sera précisé dans le théorème~\ref{thPrincMin}.

\begin{thprinc} Il existe un schéma $\agKnm$ qui est canonique, projectif de type fini sur $\Spec(\Z[1/n])$, et qui contient $\agKn$ comme ouvert dense. Il est muni d'un faisceau ample qui étend $\otimes_{i=0}^s\: \omega_i$ et l'on a
$$\agKnm \: = \: \mathrm{Proj} \left( \bigoplus_{k\in\mathbb{N}} \: \mathrm{H}^0\left( \agKn \: , \: \bigotimes_{i=0}^s \omega_i^k \right) \right)$$
Le schéma $\agKnm$ est stratifié par des variétés de Siegel de genre $\leq g$ et de niveau parahorique variable.
\end{thprinc}

Remarquons qu'utiliser uniquement $\omega_0$ conduirait à un schéma dans lequel $\agKn$ ne se plonge pas ! 
L'ingrédient principal de la démonstration du théorème est la construction de compactifications toroïdales de $\agK$~\cite{Compact@Stroh}.

Remarquons également que l'existence de la compactification minimale est bien connue pour $g=1$~\cite{Elliptique@DeligneRapoport}. Notre démarche reste valable dans ce cas, mais il faut modifier l'énoncé du théorème principal car le principe de Koecher n'est pas vérifié.

Nos résultats s'étendent à un cadre plus général que celui des variétés de Siegel. En effet, on peut remplacer $\agKn$ par une variété de Shimura P.E.L. associée à un groupe réductif $G$ sur $\Q$ non ramifié en $p$, et à un sous-groupe compact ouvert net de $G(\mathbb{A}_f)$, dont la composante en $p$ est un sous-groupe parahorique de $G(\Q_p)$ inclus dans un sous-groupe  hyper-spécial. On pourra par exemple consulter~\cite{Compact@Lan} pour les cas de bonne réduction, où la composante en $p$ est hyper-spéciale.

\subsection*{Plan de l'article}

Dans la première partie, nous expliquons comment déduire l'existence de compactifications toroïdales de $\agKn$  des résultats de~\cite{Compact@Stroh}.

Dans la seconde partie, nous étudions les formes modulaires pour $\agKn$.
Nous définissons leur développement de Fourier-Jacobi, et démontrons un principe de Koecher. Pour cela, nous utilisons le caractère réduit de $\agKn\times \Spec(\F_p)$~\cite{FlatSympl@Gortz}, la densité de son lieu ordinaire~\cite{VarAb@GenestierNgo} et une description de ses composantes irréductibles~\cite{Siegel@Yu}.

Dans la troisième partie, nous construisons le schéma $\agKnm$. Nous utilisons de manière essentielle la normalité de $\agKn$~\cite[prop.~3.1.7.1]{Compact@Stroh}.

\section{Rappels sur les compactifications toroïdales de~$\agK$}

Soient $1\leq s \leq g$ deux entiers, $p$ un nombre premier, $n\geq 3$ un entier non divisible par~$p$ et $\mathcal{D}=\{d_1<d_2<\cdots<d_s\}$ un sous-ensemble ordonné de $\{1,\cdots,g\}$.
Notons $\ag$ le champ algébrique sur $\Spec(\Z)$ qui paramètre les schémas abéliens principalement polarisés de genre~$g$. Soit $V=\oplus_{j=1}^{2g}\: \Z \cdot x_j$ un~$\Z$-module libre de rang~$2g$ muni d'une base $(x_j)$ et de la forme symplectique de matrice 
$$J=\left( \begin{array}{cc} 0 & J' \\
 -J' & 0           \end{array} \right) \: ,$$
où~$J'$ désigne la matrice  anti-diagonale dont  les coefficients non nuls sont égaux à~$1$.
Notons $\agn \rightarrow \Spec\left(\Z[1/n]\right)$ l'espace de modules des schémas abéliens principalement polarisés de genre $g$ munis d'une structure de niveau \textit{principale} en $n$. Il paramètre  les schémas abéliens~$A$  munis d'un isomorphisme
$V \: / \: nV \isolong A[n]$
qui respecte les formes symplectiques à un scalaire près. D'après~\cite{GIT@Mumford} et~\cite{Jacobi@Serre}, $\agn$ est un schéma quasi-projectif sur $\Spec(\Z[1/n])$. Notons
$\agK\rightarrow \Spec(\Z)$
le champ de modules des schémas abéliens principalement polarisés de genre $g$ munis d'une structure de niveau \textit{parahorique} de type $\mathcal{D}$ en $p$. Ce champ paramètre les schémas abéliens $A$
munis d'un drapeau de sous-groupes finis, plats et totalement isotropes 
$$ H_1 \: \subset \: H_2 \: \subset \: \cdots \: \subset \: H_s \: \subset \: A[p] $$
tel que $H_i$ soit de rang $p^{d_i}$ pour $0\leq i\leq s$. Posons enfin
$$\agKn \: = \:  \agK \times_{\ag} \agn\: .$$
Ce schéma quasi-projectif sur $\Spec(\Z[1/n])$ est l'objet central de notre article. Il paramètre les schémas abéliens principalement polarisés de genre $g$ munis d'une structure de niveau principale en $n$ et d'une structure de niveau parahorique de type $\mathcal{D}$ en $p$. Il n'est pas propre sur $\Spec(\Z[1/n])$~;~on en décrit à présent des compactifications toroïdales qui dépendent du choix d'une décomposition polyédrale d'un complexe conique.

\subsection{Données combinatoires}

Pour tout sous-espace totalement isotrope $V'$ facteur direct  de $V$, notons $V'^\perp$ son orthogonal et $C(V/V'^\perp)$ le cône des formes quadratiques semi-définies positives à radical rationnel sur $(V/V'^\perp)\otimes \R$. Notons $\mathfrak{C}$ l'ensemble des sous-espaces totalement isotropes facteurs directs de $V$, et $\mathcal{C}$  le quotient  de l'union disjointe
$$\coprod_{V'\in \mathfrak{C}} C(V/V'^\perp)$$
par la relation d'équivalence engendrée par les inclusions $C(V/V''^\perp)\subset C(V/V'^\perp)$ pour tous sous-espaces totalement isotropes facteurs directs $V''\subset V'$. C'est un \textit{complexe conique} muni d'une action de~$\gsp(\Z)$.  Si l'on pose
$$X\: = \: \bigoplus_{j=1}^g\: \Z \cdot x_j\: ,$$
il existe une injection de $C(X)$ dans $\mathcal{C}$.
Soit $\Sigma$ une décomposition polyédrale $\mathrm{GL}(X)$-admissible de $C(X)$ comme définie dans~\cite[déf.~IV.2.2]{Deg@FaltingsChai}. Son orbite dans $\mathcal{C}$ sous l'action de~$\gsp(\Z)$  définit une décomposition polyédrale $\gsp(\Z)$-admissible $\mathfrak{S}$ de $\mathcal{C}$ comme dans~\cite[déf.~3.2.3.1]{Compact@Stroh}.
D'après~\cite[th.~IV.5.7 et IV.6.7]{Deg@FaltingsChai} et~\cite[th.~3.2.7.1]{Compact@Stroh}, il existe des compactifications
$$\ag \hookrightarrow \agb \: ,\: \: \agn \hookrightarrow \agnb \: \: \mathrm{et} \: \: \agK \hookrightarrow \agKb$$
associées à $\Sigma$ et à $\mathfrak{S}$. Posons 
$$\agKnb \: = \: \agKb \times_{\agb} \agnb \: .$$
C'est une compactification de $\agKn$ sur $\Spec(\Z[1/n])$.

\subsection{Stratification de $\agKnb$} Dans ce paragraphe, on décrit deux stratifications naturelles de $\agKnb$. Soit $\mathbb{V}_\mathcal{D}^\bullet$ le drapeau de $V$ tel que
$$\mathbb{V}^i_\mathcal{D}  \:= \: \bigoplus_{j=1}^{d_i}\: \Z \cdot x_j \oplus \bigoplus_{j=d_i+1}^{2g}p\: \Z \cdot x_j$$
$$\mathrm{et} \:\:\: \: \: \mathbb{V}^{s+i}_\mathcal{D}\: = \: \bigoplus_{j=1}^{2g-d_{s+1-i}} \:  \Z\cdot x_j \oplus \bigoplus_{j=2g - d_{s+1-i}+1}^{2g} p\:\Z$$
pour $1\leq i\leq s$.
Posons
$$\Gamma_{0}= \mathrm{Stab}_{\gsp(\Z)} \left(\mathbb{V}_\mathcal{D}^\bullet\right)\quad\mathrm{et} \quad  \Gamma_n = \mathrm{Ker}\left(\gsp(\Z)\rightarrow \gsp(\Z/n\Z)\right)\: .$$
D'après~\cite[th.~IV.5.7 et IV.6.7]{Deg@FaltingsChai} et~\cite[th.~3.2.7.1]{Compact@Stroh},
les champs algébriques
$\agb$, $\agnb$ et $\agKb$
sont munis d'une stratification paramétrée par $\mathfrak{S}/\gsp(\Z)$,  $\mathfrak{S}/\Gamma_n$ et $\mathfrak{S}/\Gamma_0$, ainsi que  d'une d'une stratification plus grossière paramétrée par $\mathfrak{C}/\gsp(\Z)$,  $\mathfrak{C}/\Gamma_n$ et $\mathfrak{C}/\Gamma_0$. Dans le reste de ce paragraphe, on décrit plus précisément ces stratifications. Soit $V'\in \mathfrak{C}$.
Dans \cite{Deg@FaltingsChai} et~\cite{Compact@Stroh} sont définis des variétés de Siegel ${\mathcal{A}}_{V'}\: ,\:\:{\mathcal{A}}_{V',\: n}\:\: \mathrm{et} \: \: {\mathcal{A}}_{V',\:0}$ (de type parahorique l'image de $\mathbb{V}_{\mathcal{D}}^\bullet$ dans $V'^\perp/V'$), 
des schémas abéliens
$${\mathcal{B}}_{V'} \rightarrow {\mathcal{A}}_{V'} \: ,\:\:{\mathcal{B}}_{V',\: n}\rightarrow {\mathcal{A}}_{V',\: n}\:\: \mathrm{et} \: \: {\mathcal{B}}_{V',\: 0}\rightarrow {\mathcal{A}}_{V',\: 0}\: ,$$
des $\Z$-modules libres
$S_{V'}\: , \: \: S_{V',\: n} \: \: \mathrm{et} \: \: S_{V',\: 0}\: ,$
des torseurs
$$\mathcal{M}_{V'}\rightarrow \mathcal{B}_{V'}\: , \: \: \mathcal{M}_{V',\: n} \rightarrow \mathcal{B}_{V',\: n} \: \:\mathrm{et} \: \: \mathcal{M}_{V',\: 0} \rightarrow \mathcal{B}_{V',\: 0} $$
sous les tores
$$E_{V'}=\Hom(S_{V'},\Gm)\: , \: E_{V',\: n}=\Hom(S_{V',\: n},\Gm) \: \: \mathrm{et} \: \: E_{V',\: 0}=\Hom(S_{V',\: 0},\Gm)\: ,$$
et des fibrés en plongements toriques localement de type fini
$$\overline{\mathcal{M}}_{V'}\rightarrow \mathcal{B}_{V'}\: , \: \: \overline{\mathcal{M}}_{V',\:n} \rightarrow \mathcal{B}_{V',\: n} \: \:\mathrm{et} \: \: \overline{\mathcal{M}}_{V',\:0} \rightarrow \mathcal{B}_{V',\:0} \: $$
qui sont munis d'une strate fermée canonique.
Notons $\Gamma_{V'}$, $\Gamma_{V',\: n}$ et $\Gamma_{V',\: 0}$ les stabilisateurs de~$V'$ dans $\gsp(\Z)$, dans $\Gamma_n$ et dans $\Gamma_0$.

Notons  $\mathfrak{S}_{V'}$ l'ensemble des $\sigma\in\mathfrak{S}$  inclus dans l'intérieur de $C(V/V'^\perp)$. \`A tout $\sigma\in \mathfrak{S}_{V'}$ on associe des fibrés en plongements toriques affines de type fini
$$\overline{\mathcal{M}}_\sigma \rightarrow \mathcal{B}_{V'}\: , \: \: \overline{\mathcal{M}}_{\sigma,\:n} \rightarrow \mathcal{B}_{V',\: n} \: \:\mathrm{et} \: \: \overline{\mathcal{M}}_{\sigma,\:0} \rightarrow \mathcal{B}_{V',\:0} $$
munis d'une stratification paramétrée par les faces de $\sigma$. Dans tous les cas, la strate paramétrée par $\sigma$ (la $\sigma$-\textit{strate}) est l'unique strate fermée. Notons  $\Gamma_{\sigma}$ et $\Gamma_{\sigma,0}$ les stabilisateurs de $\sigma$ dans $\gsp(\Z)$ et dans $\Gamma_0$~;~ce sont des groupes finis.

\begin{remarque2} Comme $n\geq 3$, le stabilisateur de $\sigma$ dans $\Gamma_n$ est trivial d'après le lemme de Serre~\cite{Jacobi@Serre}.
\end{remarque2}

D'après~\cite[th.~IV.5.7 et~6.7]{Deg@FaltingsChai} et~\cite[th~3.2.7.1]{Compact@Stroh}, il existe des isomorphismes
\begin{equation} \label{isoFCStr}
\widehat{\overline{\mathcal{M}}}_{V'}\: / \: \Gamma_{V'} \isolong \widehat{\overline{\mathcal{A}}}^{V'}_g \: , \: \: \widehat{\overline{\mathcal{M}}}_{V',\: n} \isolong \widehat{\overline{\mathcal{A}}}^{V'}_{g,n} \: / \: \Gamma_{V',\: n} \: \: , \: \: \widehat{\overline{\mathcal{M}}}_{V',\:o} \: / \: \Gamma_{V',\:0} \isolong  \widehat{\overline{\mathcal{A}}}^{V'}_{g,o}
\end{equation}
$$\widehat{\overline{\mathcal{M}}}_{\sigma} \: / \: \Gamma_\sigma \isolong \widehat{\overline{\mathcal{A}}}^\sigma_g \: , \: \: \widehat{\overline{\mathcal{M}}}_{\sigma,\:n} \isolong \widehat{\overline{\mathcal{A}}}^\sigma_{g,n} \: \: \mathrm{et} \: \: \widehat{\overline{\mathcal{M}}}_{\sigma,\:o} \: / \: \Gamma_{\sigma,\:0} \isolong  \widehat{\overline{\mathcal{A}}}^\sigma_{g,0}$$
entre complétés formels le long des $V'$-strates et des $\sigma$-strates.
Considérons les produits fibrés
$${\mathcal{A}}_{V',\:n,0}\:  = \: {\mathcal{A}}_{V',\:0} \times_{{\mathcal{A}}_{V'}} {\mathcal{A}}_{V',\:n} \:\: , \: \: \:{\mathcal{B}}_{V',\:n,0}\:  = \: {\mathcal{B}}_{V',\:0} \times_{{\mathcal{B}}_{V'}} {\mathcal{B}}_{V',\:n} \: , $$
$$\overline{\mathcal{M}}_{V',\:n,0} \:= \:\overline{\mathcal{M}}_{V',\:0} \times_{\overline{\mathcal{M}}_{V'}} \overline{\mathcal{M}}_{V',\:n} \:\: , \: \: \overline{\mathcal{M}}_{\sigma,\:n,0} \:= \:\overline{\mathcal{M}}_{\sigma,\:0} \times_{\overline{\mathcal{M}}_{\sigma}} \overline{\mathcal{M}}_{\sigma,\:n} \: ,$$
la somme amalgamée
$$S_{V',\:n,0} \:= \:S_{V',\:0} \:\bigoplus_{S_{V'}} \:S_{V',\:n} \: ,$$
le groupe $\Gamma_{V',\:n,0}=\Gamma_{V',\: 0} \cap \Gamma_{V',\: n}$ et le groupe $\Gamma_{n,0}=\Gamma_{n}\cap\Gamma_{0}$. Les schémas $\overline{\mathcal{M}}_{V', \: n,0}$ et $\overline{\mathcal{M}}_{\sigma,n,0}$ sont des fibrés en plongements toriques sur  ${\mathcal{B}}_{V',\:n,0}$, équivariants sous le tore
$$E_{V',\:n,0} = \Hom(S_{V',\:n,0} ,\Gm)\: .$$
La proposition suivante résulte immédiatement des isomorphismes~(\ref{isoFCStr}).

\begin{proposition2} \label{propStratif} Le champ algébrique $\agKnb$ est muni d'une stratification paramétrée par $\mathfrak{C}/\Gamma_{n,0}$. Son complété formel le long de la $V'$-strate est isomorphe au quotient
$$\widehat{\overline{\mathcal{M}}}_{V',\: n,0} \: / \: \Gamma_{V',\: n,0}$$
du complété formel de $\overline{\mathcal{M}}_{V',\: n,0}$ le long de la strate fermée canonique. Le champ $\agKnb$ est également muni d'une stratification plus fine paramétrée par $\mathfrak{S}/\Gamma_{n,0}$. Son complété formel le long de la $\sigma$-strate est isomorphe au complété formel
$$\widehat{\overline{\mathcal{M}}}_{\sigma,n,0}$$  de ${\overline{\mathcal{M}}}_{\sigma,n,0}$ le long de la $\sigma$-strate.
\end{proposition2}

Le corollaire suivant résulte du fait  que $\agn$ et ${\overline{\mathcal{M}}}_{\sigma,n,0}$ sont des schémas pour tout $\sigma\in \mathfrak{S}$.

\begin{corollaire2} Le champ algébrique $\agKnb$ est un espace algébrique.
\end{corollaire2}

\subsection{Extension des schémas abéliens} Notons $G_0$ le schéma abélien universel sur $\agKn$ et $H_\bullet$ le drapeau universel de $G_0[p]$. Posons $G_i = G_0/ H_i$ pour $1\leq i\leq s$. Il existe une chaîne d'isogénies
$$G_0 \longrightarrow G_1 \longrightarrow \cdots \longrightarrow G_s$$
entre schémas abéliens sur $\agKn$. Le théorème~\cite[3.2.7.1]{Compact@Stroh} et les résultats de~\cite[II.2]{Deg@FaltingsChai} montrent que $G_i$ s'étend en un schéma semi-abélien
$$G_i \longrightarrow \agKnb$$
pour $0\leq i\leq s$, et que la chaîne d'isogénies s'étend en une chaîne de morphismes entre schémas semi-abéliens. Pour tout $0\leq i \leq s$, on note $e_i$ la section neutre de $G_i$ et l'on pose
$$\Omega_i=e_i^* \:\Omega_{G_i/\agKnb} \: \: \: \mathrm{et} \: \: \omega_i=\mathrm{det}(\:\Omega_i\:)\: .$$
Ces derniers sont des faisceaux localement libres sur $\agKnb$ de rang respectifs $g$ et $1$. On appelle~$\omega_i$  le $i$-ème \textit{fibré de Hodge}.  Il existe une suite de morphismes
$$\omega_s\: \longrightarrow \:\omega_{s-1} \:\longrightarrow \:\cdots \:\longrightarrow \:\omega_0$$
définis sur $\agKn$, qui induisent des isomorphismes après changement de base par $\Spec(\Z[1/n p])$.

Pour tous $V'\in\mathfrak{C}$ et $0\leq i \leq s$, il existe un schéma semi-abélien
$\tilde{G}_{V',\: i} \longrightarrow \mathcal{B}_{V',\:n,0}$ qui est globalement extension
$$0 \longrightarrow T_{V',\: i} \longrightarrow \tilde{G}_{V',\: i} \longrightarrow A_{V',\: i} \longrightarrow 0$$
d'un schéma abélien par un tore~\cite[par.~1.4.5 et~1.4.9]{Compact@Stroh}. Le schéma abélien $A_{V',\: i}$ provient d'un schéma abélien sur $\mathcal{A}_{V',\: n,0}$. Le groupe des caractères $X_{V',\: i}$ du tore $T_{V',\: i}$ est égal à l'image de $\mathbb{V}^{2s-i}_\mathcal{D}$ dans $V/V'^\perp$~\cite[par.~3.2.1]{Compact@Stroh}. On a par exemple $X_{V',0}=V/V'^\perp$. Pour tout $\sigma\in\mathfrak{S}_{V'}$, on note encore $\tilde{G}_{V',\: i}$ le schéma semi-abélien qui est l'image réciproque de $\tilde{G}_{V',\:i}$ par le morphisme
$$\widehat{\overline{\mathcal{M}}}_{\sigma,\:n,0}\longrightarrow  \mathcal{B}_{V',\:n,0}\: .$$
D'après la remarque qui suit~\cite[th.~3.1.8.2]{Compact@Stroh}, l'extension de Raynaud de
$$G_i \:\longrightarrow \: \widehat{\overline{\mathcal{M}}}_{\sigma,n,0}$$
est $\tilde{G}_{V', \: i}$. Il existe donc sur $\widehat{\overline{\mathcal{M}}}_{\sigma,n,0}$ des isomorphismes
\begin{eqnarray}\label{equationOmega}
\omega_i &\isolong & \mathrm{det} \left( \Omega_{\tilde{G}_{V',\:i}}\right) \\
\nonumber & \isolong & \mathrm{det}\left(X_{V',\:i} \right)\otimes \omega_{A_{V',\:i}}\: .
\end{eqnarray}

\section{Formes modulaires de Siegel}

Dans cette partie, on définit les formes modulaires de genre $g$, de niveau principal en $n$, et de niveau parahorique de type $\mathcal{D}$ en $p$. Le poids d'une telle forme est une famille d'entiers relatifs $(k_0,\cdots,k_s)$ qui code le faisceau inversible dont la forme est une section globale. On définit le développement de Fourier-Jacobi d'une forme modulaire, qui dépend \textit{a priori} d'un élément de $\mathfrak{S}/\Gamma_{n,0}$. On montre qu'il ne dépend en fait  que de l'élément correspondant dans l'ensemble $\mathfrak{C}/\Gamma_{n,0}$, qui est par définition l'ensemble des \textit{pointes} ou \textit{cusps}. On montre un principe de Koecher, valable si $g\geq 2$, ainsi qu'un principe de développement de Fourier-Jacobi.

Dans cette partie, on fixe un élément $\underline{k}=(k_0,\cdots,k_s)\in\Z^{s+1}\:$, un $\Z[1/n]$-module $M$, et l'on note $$\omega^{\underline{k}} \:= \:\bigotimes_{i=0}^s \: \omega_i^{\otimes k_i}$$
qui est un faisceau inversible sur $\agKnb$.

\begin{definition} Une forme modulaire de Siegel de genre $g$, de poids $\underline{k}$, de niveau principal en $n$, de niveau parahorique de type $\mathcal{D}$ en $p$, et à valeurs dans $M$ est un élément du goupe
$$\mathrm{H}^0(\agKnb \: , \: \omega^{\underline{k}}\otimes M)\: .$$
\end{definition}

\begin{remarque} La famille $\underline{k}=(k_0,\cdots,k_s)$ induit une famille de caractères $$\left(\mathrm{det}^{k_0},\cdots,\mathrm{det}^{k_0}\right)$$ du groupe algébrique $\mathrm{GL}_g$. On pourrait remplacer le choix de $\underline{k}$ par celui d'une famille de $s+1$ représentations algébriques de $\mathrm{GL}_g$, et l'on obtiendrait des formes à valeurs vectorielles~;~les résultats de cette partie s'adapteraient sans peine.
\end{remarque}

Si $M$ est un $\Z[1/pn]$-module, les faisceaux $\omega_i^{\otimes k_i}\otimes M$ sont tous isomorphes. La donnée de $\underline{k}$ est donc équivalente à celle de $\sum_i k_i$.

\begin{proposition} \label{propH0indepSigma} Le groupe $\mathrm{H}^0(\agKnb \: , \: \omega^{\underline{k}}\otimes M)$ est indépendant du choix de la décomposition polyédrale $\Sigma$.
\end{proposition}

\begin{demovide} Soit  $\Sigma'$ une autre décomposition polyédrale $\mathrm{GL}(X)$-admissibles de $C(X)$. Notons
$$\agKnb '$$
la compactification toroïdale qui lui est associée. Deux décompositions admettent un raffinement commun~;~on peut donc supposer que $\Sigma'$ raffine $\Sigma$. Il existe alors un morphisme
$$\phi \: : \: \agKnb' \longrightarrow \agKnb$$
qui est propre, surjectif et birationnel. Localement pour la topologie étale, $\phi$ s'identifie à un morphisme entre plongements toriques équivariants sous le même tore. En utilisant les résultats de~\cite[I.3]{Toro@KKMS}, on voit que
$$\phi_*\: \mathcal{O}_{\agKnb'}\: =\: \mathcal{O}_{\agKnb} $$
et que pour tout $j>0$, on a
$$R^j \: \phi_* \: \mathcal{O}_{\agKnb'} \: =\:0\: .$$
Notons $\omega'_i$ le $i$-ème fibré de Hodge sur $\agKnb'$.
Comme $\omega'_i=\phi^* \omega_{i}$, la formule de projection et la suite spectrale de Leray montrent qu'on a bien $$\mathrm{H}^0\left(\agKnb \: ,\: \omega^{\underline{k}}\otimes M\right)\: = \: \mathrm{H}^0\left(\agKnb' \: ,\: {\omega'\:}^{\underline{k}}\otimes M\right) \: . \: \:\: \qedsymbol$$
\end{demovide}

\subsection{Développement de Fourier-Jacobi}

Le $E_{V',\:n,0}\:$-torseur $\mathcal{M}_{V',\:n,0}$ est affine sur $\mathcal{B}_{V',\:n,0}$ pour tout $V'\in\mathfrak{C}$. Il existe donc un faisceau $\mathcal{A}$  en algèbres sur $\mathcal{O}_{\mathcal{B}_{V',\:n,0}}$ tel que
$$\mathcal{M}_{V',\:n,0} \:=\: \underline{\Spec}(\mathcal{A})$$
Le tore $E_{V',\:n,0}$ agit sur $\mathcal{A}$ et l'on peut décomposer ce faisceau en  algèbres selon le groupe des caractères $S_{V',\:n,0}$ de $E_{V',\:n,0}$. On obtient
$$\mathcal{A}\: = \bigoplus_{\lambda \in S_{V',\:n,0}} \mathcal{L}(\lambda)$$
où $\mathcal{L}(\lambda)$ est le sous $\: \mathcal{O}_{\mathcal{B}_{V',\:n,0}}$-module de $\mathcal{A}$ sur lequel $E_{V',\:n,0}$ agit par $\lambda$. Comme $\mathcal{M}_{V',\:n,0}$ est un $E_{V',\:n,0}\:$-torseur, $\mathcal{L}(\lambda)$ est un faisceau inversible sur $\mathcal{B}_{V',\:n,0}$. Pour tout $\sigma\in\mathfrak{S}_{V'}$, si l'on pose $$\sigma^{\vee}\: =\: \left\{\:\lambda \in  \mathrm{Sym}^2(V/V'^\perp \otimes \R) \:\: \mathrm{tel \: que} \:\: \:  b(\lambda) \geq 0\: \: \forall\: b\in\sigma\: \right\}\: ,$$
on a par définition
$$ \overline{\mathcal{M}}_{\sigma,n,0} \:=\: \underline{\Spec}\left( \bigoplus_{\lambda \in S_{V',n,0}\cap \: \sigma^\vee} \mathcal{L}(\lambda)\right)\: .$$
Posons
$$\omega_{\mathrm{et}}^{\underline{k}}\: = \: \bigotimes_{i=0}^s\: \left(\mathrm{det}(X_{V',\:i})\right)^{k_i}\quad \mathrm{et}\quad \omega_{\mathrm{ab}}^{\underline{k}}\: =\: \bigotimes_{i=0}^s \: \omega_{A_{V',\:i}}^{k_i}\: .$$
Dans la définition suivante, on utilise la proposition~\ref{propStratif} et le calcul (\ref{equationOmega})~de $\omega_i$ au voisinage d'une strate. Soit $f\in \mathrm{H}^0(\agKnb \: , \: \omega^{\underline{k}}\otimes M)$.

\begin{definition2} \label{defFJ} Le développement de Fourier-Jacobi  de $f$ en $\sigma\in\mathfrak{S}_{V'}$ est la série formelle
$$\mathrm{FJ}_{\sigma}(f) \:\:= \sum_{\lambda \in S_{V',n,0} \: \cap \: \sigma^{\vee}}  a_{\lambda,\sigma} (f)  \:\: \in \prod_{\lambda \in S_{V',n,0} \: \cap \: \sigma^{\vee}} \mathrm{H}^0\left( \mathcal{B}_{V',\:n,0}\: ,\: \mathcal{L}(\lambda)\otimes \omega_{\mathrm{et}}^{\underline{k}} \otimes \omega_{\mathrm{ab}}^{\underline{k}} \otimes M \right)$$
obtenue en évaluant $f$ sur le complété formel de $\overline{\mathcal{M}}_{\sigma,n,0}$ le long de sa $\sigma$-strate.
\end{definition2}

Posons
$\sigma^{>}\: =\: \left\{\: \lambda \in  \mathrm{Sym}^2(V/V'^\perp \otimes \R)  \:\: | \:\:  b(\lambda) > 0 \: \: \forall \: b\in\sigma \: \right\}\: .$
L'idéal de la $\sigma$-strate de $\overline{\mathcal{M}}_{\sigma,n,0}$ est engendré par les $\mathrm{H}^0(\mathcal{B}_{V',\:n,0},\mathcal{L}(\lambda))$ où $\lambda\in S_{V',\:n,0}\: \cap \:\sigma^>$. On en déduit que presque tous les $a_{\lambda,\sigma}(f)$ sont nuls pour $\lambda \in \sigma^\vee - \sigma^>$.

\begin{proposition2} Le développement de Fourier-Jacobi de $f$ en $\sigma\in\mathfrak{S}_{V'}$ ne dépend que de $f$ et de $V'$.
\end{proposition2}

\begin{demovide} Soient $\sigma,\: \tau \in \mathfrak{S}_{V'}$. Quitte à considérer une suite de cônes reliant $\tau$ et $\sigma$, on peut supposer que $\tau$ est une face de $\sigma$. Il existe alors une immersion ouverte
$$\overline{\mathcal{M}}_{\tau,n,0} \:\hookrightarrow \:\overline{\mathcal{M}}_{\sigma,n,0}\: .$$
Comme les schémas semi-abéliens $G_{V',\:i}$ se correspondent \textit{via} ce plongement, on voit que
$$\mathrm{FJ}_{\sigma}(f)\: =\: \mathrm{FJ}_{\tau}(f)\: \in\:  \mathrm{H}^0\left(\:\widehat{\overline{\mathcal{M}}}_{\tau,n,0} \: ,\: \omega^{\underline{k}} \right)\: . \:\:\: \qedsymbol$$
\end{demovide}

La proposition précédente permet de définir le développement de Fourier-Jacobi $$\mathrm{FJ}_{V'}(f)\:=\:\sum_\lambda a_{\lambda,V'}(f)$$ de $f$ en la {pointe} $V'\in \mathfrak{C}$.
Le corollaire suivant résume quelques propriétés élémentaires du développement de Fourier-Jacobi. On y utilise le fait que $\cap_{\sigma \in \mathfrak{S}_{V'}} \: \sigma^\vee = C(V/V'^\perp)^\vee$ et que le groupe $\Gamma_{V',\: n,0}$  agit sur l'algèbre
$$\prod_{\lambda \in S_{V', n,0} } \mathrm{H}^0\left(\: \mathcal{B}_{V',\: n,0}\: ,\: \mathcal{L}(\lambda)\otimes \omega_{\mathrm{et}}^{\underline{k}} \otimes \omega_{\mathrm{ab}}^{\underline{k}}\otimes M\: \right)\: .$$
\begin{corollaire2} \label{coroFJ} Le développement de Fourier-Jacobi de $f$ en $V'$ est indépendant de $\Sigma$, et est invariant sous l'action de $\Gamma_{V',\: n,0}$. L'application $\lambda\mapsto a_{\lambda,V'}(f)$ est à support dans $$S_{V',\: n,0} \: \cap\: C(V/V'^\perp)^\vee\: .$$
On a $F_{\gamma \cdot V'}\: (f) = \gamma \cdot F_{V'}(f)$ pour tout $\gamma\in \Gamma_{n,0}$.
\end{corollaire2}

\begin{remarque2} Si $V'\in \mathfrak{C}$ est un sous-espace totalement isotrope maximal, le schéma $\mathcal{B}_{V',\: n,0}$ est égal à $\Spec(\Z[1/n])$ et le développement de Fourier-Jacobi de $f$ est un élément de
$$\Z[1/n]\left[\left[\: S_{V',\: n,0} \otimes \: \omega_{\mathrm{et}}^{\underline{k}} \otimes M\: \right]\right]\: ,$$
qu'on appelle habituellement  $q$-\textit{développement} de $f$.
\end{remarque2}

\subsection{Principe de Koecher} Dans ce paragraphe, on démontre un principe de Koecher en utilisant la méthode de~\cite[V.1]{Deg@FaltingsChai} et la proposition suivante, qui résulte de~\cite{VarAb@GenestierNgo} et de~\cite{Siegel@Yu}.

\begin{proposition2} \label{propDens} Toute composante irréductible de $\agKn\times \Spec(\F_p)$ rencontre le bord de $\agKnb\times \Spec(\F_p)$.
\end{proposition2}

\begin{demo} Il suffit de démontrer l'énoncé analogue où $\F_p$ est remplacé par une clôture algébrique $\overline{\F}_p$. D'après le théorème principal de~\cite{VarAb@GenestierNgo}, le lieu ordinaire $\agKn^{\mathrm{ord}}$ de $\agKn\times \Spec(\overline{\F}_p)$ est dense dans $\agKn\times \Spec(\overline{\F}_p)$. On en déduit que les composantes irréductibles de $\agKn\times \Spec(\overline{\F}_p)$ sont les adhérences des composantes irréductibles de $\agKn^{\mathrm{ord}}$. D'après~\cite[3.2 et~3.3]{Siegel@Yu}, les composantes irréductibles de $\agKn^{\mathrm{ord}}$ sont ses composantes connexes. De plus, on dispose d'une description explicite de ces composantes connexes~:~ce sont les fibres du morphisme
$$\agKn^{\mathrm{ord}}\: \longrightarrow \: (\Z/n\Z)^*  \times \left\lbrace \: (d'_1<d'_2<\cdots<d'_s) \: \: | \: \: 0\leq  d'_i \leq d_i \:\: \forall \: i \leq s \: \right\rbrace$$
qui à un objet $\left(A,c:(\Z/n\Z)^{2g}\iso A[n],H_\bullet\right)$ de $\agKn^{\mathrm{ord}}$ associe le facteur de similitude de $c$ et la famille formée des rangs multiplicatifs $d'_i$ des $H_i$. On déduit de la proposition~\ref{propStratif} que chaque composante connexe de $\agKn^{\mathrm{ord}}$ rencontre le bord de $\agKnb\times \Spec(\overline{\F}_p)$~\cite[1.2.6]{Compact@Stroh}.
\end{demo}

\begin{remarque2} \label{remDimMax} En fait, on a démontré que toute composante irréductible de $\agKn\times \Spec(\F_p)$ rencontre une $V'$-strate pour un certain $V'\in\mathfrak{C}$ maximal.
\end{remarque2}

\begin{corollaire2}\label{coroInj} Si $FJ_{V'}(f)=0$ pour tout $V'\in\mathfrak{C}$, on a $f=0$.
\end{corollaire2}

\begin{demo} Sur un schéma quasi-cohérent, la formation de la cohomologie cohérente commute aux limites inductives de faisceaux. On peut donc supposer que $M$ est de type fini sur $\Z$, puis que $M$ est égal à $\Z/p\Z$, à $\Z[1/n]$ ou à $\Z/l\Z$, où $l$ est un nombre premier qui ne divise pas $n p$. En particulier, $M$ est un anneau et l'on a
$$\mathrm{H}^0\left(\agKnb \:,\:\omega^{\underline{k}}\otimes M\right)\: = \: \mathrm{H}^0\left(\agKnb\times \Spec(M) \:,\:\omega^{\underline{k}}\right)\: .$$
Si $M=\Z[1/n]$ ou $\Z/l\Z$, pour tout $\sigma\in \mathfrak{S}$ l'image du morphisme  $$\widehat{\overline{\mathcal{M}}}_{\sigma,n,0} \times \Spec(M)\:  \longrightarrow \: \agKnb\times \Spec(M)$$
donné par la proposition~\ref{propStratif} a une intersection non vide avec  toutes les composantes irréductibles de $\agKnb\times \Spec(M)$. Pour conclure, il suffit de remarquer que pour tout schéma réduit $X$, pour tout ouvert dense $U$ et pour tout faiseau  $\mathcal{F}$ localement libre sur $X$, la restriction
$$\mathrm{H}^0(X,\mathcal{F})\longrightarrow \mathrm{H}^0(U,\mathcal{F})$$
est injective. Si $M=\Z/p\Z$, un raisonnement analogue s'applique. On utilise la proposition~\ref{propDens} et le caractère réduit de~$\agKnb\times \Spec(M)$ \cite{FlatSympl@Gortz}, en tenant compte des développements de Fourier-Jacobi en tous les $V'\in \mathfrak{C}$.
\end{demo}

\begin{remarque2} Soit $R$ un système de représentants  des orbites de sous-espaces totalement isotropes maximaux facteurs directs de $V$ sous l'action de $\Gamma_{n,0}$. D'après la remarque~\ref{remDimMax}, l'application $f\mapsto\lbrace \mathrm{FJ}_{V'}(f) \rbrace_{V'\in R}$ est injective. Si l'on se restreint aux formes modulaires à coefficients dans un module  sans $p$-torsion, il résulte de la démonstration du corollaire~\ref{coroInj} que $f\mapsto \mathrm{FJ}_{V'}(f)$ est injective pour tout  $V'\in\mathfrak{C}$.
\end{remarque2}

On déduit du corollaire~\ref{coroInj} un principe de développement de Fourier-Jacobi.

\begin{corollaire2} \label{coroPrincDev} Soient $M\subset M'$  deux $\Z$-modules et $f\in \mathrm{H}^0(\agKnb\:\: ,\: \omega^{\underline{k}}\otimes M')$. Si $FJ_{V'}(f)$ est à coefficients dans $M$ pour tout  $V'\in\mathfrak{C}$, on a $$f\in\mathrm{H}^0(\agKnb\: ,\: \omega^{\underline{k}}\otimes M)\: .$$
\end{corollaire2}

On dispose également d'un principe de Koecher.

\begin{proposition2} \label{propPrincKoech} Si $g\geq 2$, l'application de restriction réalise un isomorphisme
$$\mathrm{H}^0\left(\agKnb \: ,\: \omega^{\underline{k}}\otimes M\right)\:\isolong \:   \mathrm{H}^0\left(\agKn\: ,\: \omega^{\underline{k}}\otimes M\right)\: .$$
\end{proposition2}

\begin{demo} On se réduit comme précédemment au cas où $M=\Z[1/n]$ ou $\Z/l\Z$, et où $l$ ne divise pas $n$. Le schéma $\agKn\times \Spec(M)$ est un ouvert dense de $\agKnb\times \Spec(M)$. Le schéma $\agKnb\times \Spec(M)$ est réduit. On en déduit que l'application de restriction est injective.
Montrons qu'elle est surjective. Soit
$f\in \mathrm{H}^0(\agKn\times \Spec(M) \: ,\: \omega^{\underline{k}})$. Pour tout $V'\in \mathfrak{C}$  maximal, on peut former un développement de Fourier-Jacobi $\mathrm{FJ}_{V'}(f)$ analogue à celui de la définition~\ref{defFJ}. Notons $\Lambda$ l'ensemble des $\lambda\in  S_{V',\: n,0}$ tels que $a_{\lambda,V'}(f)\neq 0$. Il est invariant sous  $\Gamma_{V',\: n,0}$ et son intersection avec $S_{V',\: n,0}\setminus\sigma^\vee$ est finie pour tout $\sigma\in \mathfrak{S}_{V'}$. On déduit de~\cite[lem.~V.1.3]{Deg@FaltingsChai} que $\Lambda \subset C(V/V'^\perp)^\vee$~:~en effet, on remarque que l'image de $\Gamma_{V',\: n,0}$ dans $\mathrm{GL}(V/V'^\perp)$ est d'indice fini, et qu'une puissance assez grande de l'endomorphisme~$T$ de~\cite[V.1.3]{Deg@FaltingsChai} est dans cette image. La section $f$ s'étend donc à $\widehat{\overline{M}}_{V',n,0}$. On en déduit que~$f$ s'étend $\widehat{\overline{M}}_{V'',n,0}$ pour tout $V''\in\mathfrak{C}$, puis que $f$ s'étend à $\agKnb$.
\end{demo}

\subsection{Coefficient constant}

Soient une forme modulaire $f\in \mathrm{H}^0(\agKnb\: ,\: \omega^{\underline{k}}\otimes M)$, un sous-espace  $V'\in\mathfrak{C}$ et un cône $\sigma\in\mathfrak{S}_{V'}$. Dans cette partie, nous montrons que la restriction de $f$ à la $\sigma$-strate est indépendante de $\sigma$ et définit une forme modulaire de genre $g-\mathrm{rang}(V')$.

\begin{proposition2} La restriction de $f$ à la $\sigma$-strate de $\agKnb$ ne dépend que de $V'$ et est égale à
$$a_{0,V'}(f)\: \in\: \mathrm{H}^0\left(\mathcal{B}_{V',\:n,0}\: ,\:  \omega_{\mathrm{et}}^{\underline{k}} \otimes \omega_{\mathrm{ab}}^{\underline{k}}\otimes M\right)$$
\end{proposition2}

\begin{demo} L'idéal de la $\sigma$-strate est engendré par les $\mathrm{H}^0(\mathcal{B}_{V',\:n,0}\:,\mathcal{L}(\lambda))$ où $\lambda\in S_{V',\:n,0} \: \cap \:  \sigma^>$, et le développement de Fourier-Jacobi de $f$ est à support dans $C(V/{V'}^\perp)^\vee$. Mais $C(V/{V'}^\perp)^\vee \: \setminus \sigma^> = \{0 \}$ car $\sigma$ est inclus dans l'intérieur de $C(V/{V'}^\perp)$.
\end{demo}

On rappelle que pour $0\leq i \leq s$, le schéma abélien $A_{V',\:i}$ descend à la variété de Siegel $\mathcal{A}_{V',\:n,0}$ de genre $g-\mathrm{rang}(V')$ et de type parahorique l'image de $\mathbb{V}_{\mathcal{D}}^\bullet$ dans $V'^\perp/V'$. D'après~\cite[prop.~1.4.6.1]{Compact@Stroh}, la flèche $\mathcal{B}_{V',\:n,0}\rightarrow \mathcal{A}_{V',\:n,0}$ admet une section canonique dont l'image s'identifie aux   invariants  de $\mathcal{B}_{V',\:n,0}$ par $\Gamma_{V',\:n,0}$. On en déduit le résultat suivant.

\begin{corollaire2}\label{coroCoeffConst} La section $a_{0,V'}(f)$ provient d'une section de $\mathrm{H}^0\left(\mathcal{A}_{V',\:n,0} \: ,\: \omega_{\mathrm{ab}}^{\underline{k}}\otimes M\right)$.
\end{corollaire2}

\section{Compactification minimale}

Construisons à présent la compactification minimale  de $\agKn$. Dans un premier temps, on se restreint au cas où $n> 2p g\sqrt{2p}$. Cette hypothèse sera utilisée dans la preuve de la  proposition~\ref{propFibreGeoAb}.

Le théorème~\cite[IX.2.1]{Pinceaux@MoretBailly} affirme que si $S$ est un schéma normal excellent, si $G\rightarrow S$ est un schéma semi-abélien et si $\omega_{G/S}$ désigne le déterminant du faisceau des formes différentielles invariantes de $G$, il existe un entier positif $m$ tel que $$\omega_{G/S}^{\otimes m}$$ soit engendré par ses sections globales sur $S$.  Le schéma $\agKnb$ est  excellent et normal d'après~\cite[th.~3.2.7.1]{Compact@Stroh}. Il existe donc un entier positif $m$ tel que  $\omega_i^{\otimes m}$ soit engendré par ses sections globales pour tout $0\leq i \leq s$.
Ainsi, il existe un morphisme propre
$$\varphi \: : \: \agKnb \: \longrightarrow \: \prod_{i=0}^s \mathbb{P}\left( \mathrm{H}^0(\agKnb\: ,\omega_i^{m})\right)\: .$$
On applique le théorème de factorisation de Stein~\ega{iii}{4.3.3,\: 4.3.4} et l'on obtient un schéma $\agKnm$ et une factorisation de $\varphi$ en
$$\agKnb\: {\buildrel{\pi}\over\longrightarrow}  \: \agKnm \:{\buildrel{\psi}\over\longrightarrow} \: \prod_{i=0}^s \mathbb{P}\left(  \mathrm{H}^0(\agKnb\: ,\omega_i^{m})\right)  \: ,$$
où $\psi$ est un morphisme fini, où les fibres de $\pi$ sont non vides et géométriquement connexes, et où $\pi_*\mathcal{O}_{\agKnb}=\mathcal{O}_{\agKnm}$. Pour tout $0\leq i \leq s$, il existe un faisceau inversible $\mathcal{L}_i$ sur  $\mathbb{P}\left(  \mathrm{H}^0(\agKnb,\omega_i^{m})\right)$ qui est très ample et tel que $\varphi^* \mathcal{L}_i \iso \omega_i^m$. L'existence du plongement de Segre~\ega{ii}{4.3} montre que $\boxtimes_{i=0}^s \mathcal{L}_i$ est très ample sur $\prod_{i=0}^s \mathbb{P}\left(  \mathrm{H}^0(\agKnb,\omega_i^{m})\right) $. Posons $\mathcal{L}=\psi^* \: \boxtimes_{i=0}^s \mathcal{L}_i$~; d'après~\ega{ii}{4.6.3}, c'est un faisceau ample sur $\agKnm$.
Il résulte de~\ega{ii}{4.6.3} que
$$\agKnm \: = \: \mathrm{Proj} \left(  \bigoplus_{k\in \mathbb{N}}\: \mathrm{H}^0 ( \:\agKnm\:,\:\mathcal{L}^k\: )\right)\: .$$
Comme $\pi^* \mathcal{L}\iso \otimes \: \omega_i^m$ et $\pi_{*} \mathcal{O}_{\agKnb} = \mathcal{O}_{\agKnm}$, on a $\mathrm{H}^0 ( \agKnm,\mathcal{L}^k )= \mathrm{H}^0 ( \agKnb,
\otimes_i \:  \omega_i^{k m} )$ pour tout $k\geq 0$. Notons $\underline{k}=(k,\cdots,k)\in\mathbb{N}^{s+1}$ pour tout $k\in \mathbb{N}$. D'après~\ega{ii}{2.4.7}, on a
$$\agKnm \: = \: \mathrm{Proj} \left(  \bigoplus_{k\in \mathbb{N}}\: \mathrm{H}^0 ( \:\agKnb\:,\: \omega^{\underline{k}}\: )\right)\: .$$
En particulier, l'algèbre graduée $\bigoplus_{k\in \mathbb{N}}\: \mathrm{H}^0 ( \:\agKnb\:,\: \omega^{\underline{k}}\: )$ est de type fini sur $\Z[1/n]$. La proposition~\ref{propH0indepSigma} montre que $\agKnm$ est indépendant du choix de $\Sigma$.

\begin{lemme} \label{lemGiDetermineChaine} Soient $g$, $m$ et $N$ trois entiers tels que $m>2 N g \sqrt{2N}$, soit $K$ un corps, et soient $G$ et $G'$ deux schémas abéliens de genre $g$ sur $\Spec(K)$ munis de polarisations $\lambda$ et $\lambda'$ et de structures de niveau principales en $m$. Si $\lambda$ est principale, il existe au plus une isogénie $f \: : \:  G\rightarrow G'$  compatible aux structures de niveau en $m$ telle que $f^t \circ \lambda' \circ f = N\lambda$, où $f^t$ désigne l'isogénie duale de $f$.
\end{lemme}

\begin{demo} Soient $f$ et  $g$ deux isogénies satisfaisant les hypothèses du lemme. Comme  $\mathrm{Hom}(G,G')$ s'identifie à un réseau de $\mathrm{End}(G)\otimes \Q$ inclus dans $\mathrm{End}(G)\otimes \Z[1/N]$, on  peut voir  $f$ et $g$ comme des $\Q$-endomorphismes de $G$ tels que $Nf$ et $Ng$ soient dans $\mathrm{End}(G)$. D'après~\cite[coro.~19.3]{Abelian@Mumford}, la $\Q$-algèbre $\mathrm{End}(G)\otimes \Q$ est semi-simple de rang $\leq 4g^2$. Notons~$\mathrm{Tr}$ la trace réduite qui lui associée et $q$ la forme quadratique sur $\mathrm{End}(G)\otimes \Q$ qui envoie~$u$ sur $\mathrm{Tr}(u^t \lambda u \lambda^{-1} )$. Cette dernière est définie positive par~\cite[th.~21.1]{Abelian@Mumford} et prend des valeurs entières sur $\mathrm{End}(G)$. Par hypothèse, on~a
$$q(f) =\mathrm{Tr}(p)\leq 4N g^2\quad , \quad q(g) \leq 4Ng^2$$
et il existe $u\in \mathrm{Hom}(G,G')$ tel que $f-g=m u$. On a donc 
$$m^2 q(u) \:=\: q(f-g) \:\leq \:q(f) + q(g) \:\leq\: 8N g^2\: .$$
Comme $m>2 N g\sqrt{2N}$ par hypothèse, on en déduit que $q(u)<1/N^2$. Or on a 
$$Nu \:\in\: \mathrm{End}(G)$$
et $q$ prend des valeurs entières sur $\mathrm{End}(G)$, donc $u=0$ et $f=g$.
\end{demo}

\begin{remarque} En posant $N=1$ puis $G=G'$ et $\lambda=\lambda'$ dans le lemme précédent, on montre que tout automorphisme de $G$ préservant une polarisation principale et une structure principale de niveau $>2g\sqrt{2}$ est égal à l'identité. Ceci est une forme faible du lemme de Serre~\cite{Jacobi@Serre}.
\end{remarque}

\begin{proposition} \label{propFibreGeoAb} Les fibres géométriques de $\pi$ au-dessus de l'image de $\agKn$ possèdent un unique point géométrique.
\end{proposition}

\begin{demo} Comme $\varphi^*\mathcal{L}_i\iso \omega_i^m$, le faisceau $\omega_i^m$ est trivialisé le long des fibres de $\pi$ pour $0\leq i \leq s$. D'après~\cite[th.~X.4.5]{Pinceaux@MoretBailly} ou~\cite[prop.~V.2.2]{Deg@FaltingsChai}, on en déduit que pour tout point géométrique $\bar{x}$ de $\pi(\agKn)$ et tout $0\leq i\leq s$, la fibre $\pi^{-1}(\bar{x})$ est incluse dans $\agKn$ et la classe d'isomorphisme du schéma abélien $G_i$ est constante sur $\pi^{-1}(\bar{x})$. Le faisceau des endomorphismes de $G_i$ est  constant sur $\pi^{-1}(\bar{x})$ , à valeurs dans un $\Z$-module libre de type fini. Cela implique que la classe d'isomorphisme du schéma abélien polarisé $G_i$ est constante sur $\pi^{-1}(\bar{x})$. On rappelle qu'on suppose $n> 2p g\sqrt{2p}$. On applique le lemme~\ref{lemGiDetermineChaine} à $G=G_0$, à $G'=G_i$, à la polarisation principale $\lambda=\lambda_0$ de $G_0$, à la polarisation canonique $\lambda'=\lambda_i$ de $G_i$, à $m=n$ et à $N=p$ pour $0\leq i \leq s$. On en déduit que la classe d'isomorphisme de la chaîne d'isogénies 
$$G_0 \: \longrightarrow  \: G_1 \: \longrightarrow \: \cdots \: \longrightarrow G_s$$
entre schémas abéliens polarisés est constante sur $\pi^{-1}(\bar{x})$. Cela montre que la fibre géométrique $\pi^{-1}(\bar{x})$ possède un unique point géométrique.
\end{demo}

\begin{corollaire} \label{coroPiIsoSurOuvert} La restriction de $\pi$ à $\agKn$ est un isomorphisme sur son image.
\end{corollaire}

\begin{demo} D'après la proposition~\ref{propFibreGeoAb} et  le \textit{Main Theorem} de Zariski~\ega{iii}{4.4.3}, la restriction de $\pi$ à $\agKn$ est finie sur son image. La proposition~\ref{propFibreGeoAb} et le caractère réduit de $\agKn$ montrent que la restriction de $\pi$ à $\agKn$ est une application birationnelle. Comme $\pi_{*} \mathcal{O}_{\agKnb} = \mathcal{O}_{\agKnm}$ et $\agKnb$ est normal, le théorème des fonctions formelles~\ega{iii}{4.1.5} montre que le schéma $\agKnm$ est normal. Toute application finie birationnelle d'image normale est un isomorphisme. C'est donc le cas du morphisme $\pi:\agKn \rightarrow \pi(\agKn)$ obtenu par restriction.
\end{demo}

\'Etudions à présent la restriction de $\pi$ aux autres strates de $\agKnb$. Comme $\pi$ est donné par l'évaluation de formes modulaires de Siegel,  le corollaire~\ref{coroCoeffConst} implique que pour tout $V'\in \mathfrak{C}$, la restriction de $\pi$ à la $V'$-strate se factorise en
$$\pi_{V'} \: : \: \mathcal{A}_{V',\:n,0} \longrightarrow \agKnm\: .$$
On montre comme précédemment que $\pi_{V'}$ est un morphisme fini birationnel sur son image, dont les fibres géométriques possèdent un unique point géométrique. Avant de prouver la normalité de $\pi_{V'}(\mathcal{A}_{V',\:n,0})$, il nous faut établir quelques résultats préliminaires.

\begin{proposition}\label{propStratifAgMin} Deux strates paramétrées par des éléments différents de $\mathfrak{C}/\Gamma_{n,0}$ ont des images disjointes par $\pi$.
\end{proposition}

\begin{demo} Soit $\bar{x}$ un point géométrique de $\agKnm$. D'après~\cite[prop.~V.2.2]{Deg@FaltingsChai}, la restriction de $G_i$ à $\pi^{-1}(\bar{x})$  pour $0\leq i\leq s$ est un schéma semi-abélien de rang torique constant. Ainsi, les rangs multiplicatifs et abéliens~\cite[1.2.6]{Compact@Stroh} de la chaîne
\begin{equation}\label{eqchaine}
G_0 \rightarrow G_1 \rightarrow \cdots \rightarrow G_s
\end{equation}
sont constants sur $\pi^{-1}(\bar{x})$. Or d'après~\cite[1.2.6]{Compact@Stroh}, ces rangs paramètrent exactement l'ensemble $\mathfrak{C}/\Gamma_{n,0}$.
\end{demo}

Le schéma $\agKnm$ est donc muni d'une stratification indexée par 
$\mathfrak{C}/\Gamma_{n,0}$. Notons $Z_{V'}$ la strate de $\agKnm$ paramétrée par $V'\in \mathfrak{C}$. Soit $\bar{x}$ un point géométrique de $Z_{V'}$. Notons encore $\bar{x}$ l'unique point géométrique de la fibre de $\pi_{V'}:\mathcal{A}_{V',\: n,0} \longrightarrow \agKnm$ en $\bar{x}$. Notons $$\widehat{\mathcal{O}}_{\agKnm\:,\: \bar{x}}$$ le complété du localisé strict $\mathcal{O}_{\agKnm,\: \bar{x}}$ de $\agKnm$ en $\bar{x}$ et $$\widehat{\mathcal{B}}_{V',\: n,0,\: \bar{x}}$$ le complété formel de $\mathcal{B}_{V',\: n,0}\times \Spec(\mathcal{O}_{\mathcal{A}_{V',n,0},\: \bar{x}})$ le long de la fibre  de $\mathcal{B}_{V',\: n,0}\rightarrow \mathcal{A}_{V',\: n,0}$ en $\bar{x}$.

\begin{proposition} \label{propNormal} Il existe un isomorphisme canonique
$$\widehat{\mathcal{O}}_{\agKnm\:,\: \bar{x}} \isolong \left( \prod_{\lambda\in S_{V',n,0}\:\cap\: C(V/V'^\perp)^\vee}  \mathrm{H}^0 \left(  \widehat{\mathcal{B}}_{V',\: n,0,\: \bar{x}} \: , \: \mathcal{L}(\lambda)  \right)   \right) ^{\Gamma_{n,0}}$$
qui induit un isomorphisme
$$\widehat{\mathcal{O}}_{Z_{V'},\: \bar{x}} \isolong \widehat{\mathcal{O}}_{ \mathcal{A}_{V',\: n,0}\:,\: \bar{x}} $$
\end{proposition}

\begin{demo} D'après~\ega{iii}{4.1.5}, l'anneau complet $\widehat{\mathcal{O}}_{\agKnm,\: \bar{x}}$ est isomorphe à l'anneau des fonctions régulières sur le complété formel de $\agKnb$ le long de $\pi^{-1}(\bar{x})$. La proposition~\ref{propStratifAgMin} montre que $\pi^{-1}(\bar{x})$ est inclus dans la $V'$-strate de $\agKnb\times \Spec(\mathcal{O}_{\agKnm,\:\bar{x}})$. L'existence du premier isomorphisme résulte  de l'isomorphisme
$$\widehat{\overline{\mathcal{M}}}_{V',\: n,0} \: / \: \Gamma_{V',\: n,0} \: \isolong \: \widehat{\overline{\mathcal{A}}}^{V'}_{g,n,0}$$
et du fait que l'anneau des fonctions régulières sur $\widehat{\overline{\mathcal{M}}}_{V',\: n,0}$ est
$$\prod_{\lambda\in S_{V', n,0}\:\cap \:C(V/V'^\perp)^\vee}  \mathrm{H}^0 \left(  \widehat{\mathcal{B}}_{V',\: n,0} \: , \: \mathcal{L}(\lambda) \right) \: .$$
La seconde assertion se déduit de la première et du corollaire~\ref{coroCoeffConst}.
\end{demo}

\begin{corollaire} \label{coroStrate} Pour tout $V'\in\mathfrak{C}$, la restriction de $\pi$ à la $V'$-strate de $\agKnb$ se factorise en un isomorphisme
$\pi_{V'} : \mathcal{A}_{V',\: n,0} \:\isolong \:Z_{V'}$.
\end{corollaire}

\begin{demo} Soit $V'\in\mathfrak{C}$. On a vu que la restriction de $\pi$ à la $V'$-strate de $\agKnb$ se factorisait en un morphisme $\pi_{V'}: \mathcal{A}_{V',\: n,0} \rightarrow \agKnm$. Par définition, le sous-schéma $Z_{V'}$ de $\agKnm$ est l'image de $\pi_{V'}$. On a montré que le morphisme $\pi_{V'}:\mathcal{A}_{V',\: n,0} \rightarrow Z_{V'}$ est fini birationnel. D'après~\cite[prop.~3.1.7.1]{Compact@Stroh}, le schéma $\mathcal{A}_{V',\: n,0}$ est normal. On déduit du second isomorphisme de la proposition~\ref{propNormal} que le schéma $Z_{V'}$ est normal. Ainsi, $\pi_{V'}$ induit un isomorphisme de $\mathcal{A}_{V',\: n,0}$ sur $Z_{V'}$.
\end{demo}

Soit $n'$  un entier divisible par $n$ mais non par $p$. Le schéma $\mathcal{A}_{g,n',0}$ est muni d'une action du groupe $\gsp(\Z/n'\Z)$, et le quotient par l'action du  sous-groupe $\Gamma_{n,0}/\Gamma_{n',0}$ s'identifie canoniquement à $\agKn
$. Comme l'éventail $\Sigma$ est équivariant sous l'action de  $\mathrm{GL}(X)$, le schéma $\mathcal{A}^*_{g,n',0}$ est muni d'une action de $\gsp(\Z/n'\Z)$ prolongeant celle sur $\mathcal{A}_{g,n',0}$.  Le quotient de $\mathcal{A}^*_{g,n',0}$ par
$\Gamma_{n,0}/\Gamma_{n',0}$ est canoniquement isomorphe à $\agKnm$.

Pour tout $n\geq 3$, il existe un entier $n'> 2p g\sqrt{2p}$ qui est divisible par $n$ mais non par $p$. En posant 
$$\agKnm \: = \: \mathcal{A}^*_{g,n',0} \: / \: \Gamma_{n,0}\: ,$$
on obtient une compactification minimale de $\agKn$ qui est indépendante du choix de $n'$.

\begin{remarque} De la même manière, on obtient une compactification minimale de l'espace de modules grossier associé à $\agK$. Sa description est moins élégante que celle de $\agKnm$ car il faut tenir compte des automorphismes du champ $\agK$.
\end{remarque}

Notons $\agnm$ la compactification minimale de $\agn$ construite dans~\cite[V.2.5]{Deg@FaltingsChai}. On rappelle que $\agKnb$ désigne la compactification toroïdale associée à $\Sigma$. Dans le théorème suivant, nous résumons les propriétés de la construction que nous venons d'effectuer. Le cas des courbes elliptiques étant bien connu~\cite{Elliptique@DeligneRapoport}, on suppose $g\geq 2$.

\begin{theoreme}\label{thPrincMin} Pour tout $n\geq 3$, il existe un schéma canonique $\agKnm$ projectif de type fini sur $\Spec(\Z[1/n])$. Il contient  $\agKn$ comme ouvert dense, est normal et indépendant de $\Sigma$.  Le morphisme d'oubli de $\agKn$ vers $\agn$ s'étend en un morphisme de $\agKnm$ vers $\agnm$. 
Pour tout $0\leq i \leq s$, le fibré de Hodge $\omega_i$ sur $\agKn$ s'étend en un faisceau inversible sur $\agKnm$. Le produit tensoriel $\otimes_{i=0}^s \: \omega_i$ est ample sur $\agKnm\:$, et il existe un isomorphisme canonique
$$\agKnm \:\isolong \:\mathrm{Proj}\left( \bigoplus_{k\in \mathbb{N}} \:\mathrm{H}^0 \left( \agKn \: ,\: \bigotimes_{i=0}^s \:\omega_i^k \right)  \right) \: .$$
Le schéma $\agKnm$ est muni d'une stratification localement fermée paramétrée par $\mathfrak{C}/\Gamma_{n,0}$;  la strate associée à $V'\in\mathfrak{C}$ est canoniquement isomorphe à $\mathcal{A}_{V',\: n,0}$.
Il existe un morphisme canonique
$$\pi \: : \: \agKnb \longrightarrow \agKnm\: .$$
L'image inverse par $\pi$ de la stratification de $\agKnm$  coïncide avec la stratification de $\agKnb$ qui est paramétrée par $\mathfrak{C}/\Gamma_{n,0}$.
\end{theoreme}

\begin{demo} L'isomorphisme entre $\agKnm$ et $\mathrm{Proj}\left( \oplus_{k\in \mathbb{N}} \:\mathrm{H}^0 \left( \agKn \: ,\: \otimes_{i=0}^s \:\omega_i^k \right)  \right)$ résulte de l'isomorphisme
$$\agKnm \:\isolong\: \mathrm{Proj}\left( \bigoplus_{k\in \mathbb{N}} \:\mathrm{H}^0 \left( \agKnb \: ,\: \bigotimes_{i=0}^s \:\omega_i^k \right)  \right)$$
et du principe de Koecher~\ref{propPrincKoech}, qui est valable si $g\geq 2$.
Il reste à prouver que $\omega_i$ s'étend en un faisceau inversible sur $\agKnm$ pour $0\leq i \leq s$. Soit $j$ l'immersion de $\agKn$ dans $\agKnm$. Comme $\agKnm$ est normal et $\agKnm\setminus \agKn$ est de codimension $\geq 2$, une éventuelle extension localement libre du faisceau $\omega_i$ est unique, et égale à $j_* \: \omega_i$. Par descente, il suffit de montrer que~$\omega_i$ s'étend au complété du localisé strict de $\agKnm$ en tout point géométrique d'une $V'$-strate. Cela résulte de~(\ref{equationOmega}) et du corollaire~\ref{coroStrate}.
\end{demo}

\bibliography{bibliographie}
\bibliographystyle{smfalpha}













\end{document}

%% file: paquets.tex
\usepackage{amsmath,a4wide}
\usepackage{amssymb}
\usepackage{xspace}
\usepackage{euscript}
\usepackage{mathrsfs} 
\usepackage[all]{xy}
\usepackage{hyperref}
\usepackage[francais]{babel}
\usepackage[T1]{fontenc}
\usepackage{textcomp}
\usepackage{smfthm}

%% file: lemme-proposition-theoreme.tex
\newtheorem{lemme}[subsection]{Lemme}		

\newtheorem{proposition}[subsection]{Proposition}	
\newtheorem{proposition2}[subsubsection]{Proposition}

\newtheorem*{thprinc}{Théorème principal}

\newtheorem{theoreme}[subsection]{Th\'eor\`eme}

\newtheorem{corollaire}[subsection]{Corollaire}
\newtheorem{corollaire2}[subsubsection]{Corollaire}

\theoremstyle{definition}
\newtheorem{definition}[subsection]{D\'efinition}
\newtheorem{definition2}[subsubsection]{D\'efinition}

\newtheorem{remarque}[subsection]{Remarque}
\newtheorem{remarque2}[subsubsection]{Remarque}

\newenvironment{demo}{\noindent{\textit{D\'emonstration. --- }}}{~\qedsymbol \vspace{4mm}}
\newenvironment{demovide}{\noindent{\textit{D\'emonstration. --- }}}{}

%% file: numerotation.tex
\renewcommand{\theenumi}{\textbf{\roman{enumi}}} 

\renewcommand\theequation{\thesection.\textbf{\alph{equation}}} 

\makeatletter
\@addtoreset{equation}{subsection}
\makeatother

\renewcommand\thesubsubsection{\textbf{\thesubsection.\arabic{subsubsection}}}
\renewcommand\thesubsection{\textbf{\thesection.\arabic{subsection}}}

%% file: macros.tex
\newcommand{\Z}{\mathbb Z}
\newcommand{\Q}{\mathbb Q}
\newcommand{\R}{\mathbb R}

\newcommand{\F}{\mathbb F}

\newcommand{\Gm}{{\mathbb G}_m}

\newcommand{\Hom}{\underline{\mathrm{Hom}}}

\newcommand{\Spec}{\mathrm{Spec}}

\newcommand{\iso}{\buildrel{\sim}\over\rightarrow}

\renewcommand{\L}{\mathrm{L}}		

\def\commutatif{\ar@{}[rd]|{\circlearrowleft}}
\def\cartesien{\ar@{}[rd]|{\square}}

\renewcommand{\lim}{{\mathrm{lim}}} 

\def\ega#1#2{[{\bf \'EGA}~{\sc #1}~#2]} 				
\def\egazéro#1#2{[{\bf \'EGA}~$0_{\textsc{#1}}$~#2]}		
